\documentclass[a4paper,reqno]{amsart}
\usepackage[margin=2.5cm,papersize={17.75cm,24.8cm}]{geometry}
\usepackage[all]{xy}           
\usepackage{amssymb}           
\usepackage{hyperref}
\usepackage{eucal}
\usepackage[dvips]{graphics}
\usepackage{graphicx}
\usepackage{exscale}
\usepackage{color}

\numberwithin{equation}{section}

\newtheorem{definition}{Definition}[section]

\newtheorem{theorem}[definition]{Theorem}
\newtheorem{proposition}[definition]{Proposition}
\newtheorem{corollary}[definition]{Corollary}
\newtheorem{remarkth}[definition]{Remark}
\newtheorem{example}[definition]{Example}
\newenvironment{remark}{\begin{remarkth}\upshape}{\hfill$\diamond$\end{remarkth}}

\renewcommand{\emph}[1]{{\bfseries\itshape{#1}}}

\newcommand{\R}{\mathbb{R}}      
\newcommand{\Z}{\mathbb{Z}}      
\newcommand{\F}{\mathbb{F}}

\newcommand{\J}{\mathbb{J}}
\newcommand{\lcf}{\lbrack\! \lbrack}
\newcommand{\rcf}{\rbrack\! \rbrack}

\newcommand{\lvec}[1]{\overleftarrow{#1}}
\newcommand{\rvec}[1]{\overrightarrow{#1}}

\makeatletter
\newcommand\prol{\@ifstar{\@proldf}{\@prolpf}}  
\def\@prolpf{\@ifnextchar[{\@prolpf@wrt}{\@prolpf@}}
\def\@prolpf@wrt[#1]#2{\@ifnextchar[{\@prolpf@wrt@at{#1}{#2}}{\@prolpf@wrt@{#1}{#2}}}
\def\@prolpf@wrt@at#1#2[#3]{\prolsymbol^{#1}_{#3}#2}
\def\@prolpf@wrt@#1#2{\prolsymbol^{#1}#2}
\def\@prolpf@#1{\@ifnextchar[{\@prolpf@at{#1}}{\@prolpf@@{#1}}}
\def\@prolpf@at#1[#2]{\prolsymbol_{#2}#1}
\def\@prolpf@@#1{\prolsymbol#1}
\def\@proldf{\@ifnextchar[{\@proldf@wrt}{\@proldf@}}
\def\@proldf@wrt[#1]#2{\@ifnextchar[{\@proldf@wrt@at{#1}{#2}}{\@proldf@wrt@{#1}{#2}}}
\def\@proldf@wrt@at#1#2[#3]{\prolsymbol^{*#1}_{#3}#2}
\def\@proldf@wrt@#1#2{\prolsymbol^{*#1}#2}
\def\@proldf@#1{\@ifnextchar[{\@proldf@at{#1}}{\@proldf@@{#1}}}
\def\@proldf@at#1[#2]{\prolsymbol^*_{#2}#1}
\def\@proldf@@#1{\prolsymbol^*#1}
\def\prolsymbol{\mathcal{T}}
\makeatother

\def\lcf{\lbrack\! \lbrack}
\def\rcf{\rbrack\! \rbrack}
\newcommand{\setdef}[2]{\{#1 \; | \; #2\}}

\newcommand{\set}[2]{\left\{\,#1\left.\vphantom{#1#2}\,\right\vert\,#2\,
                \right\}}

\def\gpd{\rightrightarrows}
\def\lcf{\lbrack\! \lbrack}
\def\rcf{\rbrack\! \rbrack}
\def\gpd{\rightrightarrows}

\begin{document}

\title[Discrete dynamics  in implicit form]{Discrete dynamics  in implicit form}

\author[D.\ Iglesias]{David Iglesias}
\address{David Iglesias:
Unidad asociada ULL-CSIC ``Geometr\'{\i}a Diferencial y Mec\'anica Geo\-m\'e\-tri\-ca". Departamento de Matem\'atica Fundamental, Facultad de
Ma\-te\-m\'a\-ti\-cas, Universidad de la Laguna, La Laguna,
Tenerife, Canary Islands, Spain} \email{diglesia@ull.es}

\author[J.\ C.\ Marrero]{Juan C.\ Marrero}
\address{Juan C.\ Marrero:
Unidad asociada ULL-CSIC ``Geometr\'{\i}a Diferencial y Mec\'anica Geo\-m\'e\-tri\-ca".
Departamento de Matem\'atica Fundamental, Facultad de
Ma\-te\-m\'a\-ti\-cas, Universidad de la Laguna, La Laguna,
Tenerife, Canary Islands, Spain} \email{jcmarrer@ull.es}

\author[D.\ Mart{\'\i}n de Diego]{D. Mart{\'\i}n de Diego}
\address{D. Mart{\'\i}n de Diego:
Instituto de Ciencias Matem\'aticas, CSIC-UAM-UC3M-UCM,
Campus de Cantoblanco, UAM,
C/Nicol\'as Cabrera, 15
 28049 Madrid, Spain}
\email{david.martin@icmat.es}

\author[E.\  Padr\'on]{Edith Padr\'on}
\address{Edith Padr\'on:
Unidad asociada ULL-CSIC ``Geometr\'{\i}a Diferencial y Mec\'a\-ni\-ca Geo\-m\'e\-tri\-ca".
Departamento de Matem\'atica Fundamental, Facultad de
Ma\-te\-m\'a\-ti\-cas, Universidad de la Laguna, La Laguna,
Tenerife, Canary Islands, Spain} \email{mepadron@ull.es}

\thanks{This work has been partially supported by MICINN (Spain)
Grants  MTM2009-13383,  MTM2010-21186-C02-01 and  MTM2009-08166-E, project "Ingenio
Mathematica" (i-MATH) No. CSD 2006-00032 (Consolider-Ingenio 2010), the project of the
Canary Islands government SOLSUB200801000238 and the European project IRSES-project ``Geomech-246981''.
D.~Iglesias wishes to thank MICINN for a ``Ram\'on y Cajal" research
contract.}
\begin{abstract}
A notion of implicit difference equation on a Lie groupoid is introduced
and an algorithm for extracting the integrable part (backward or/and forward) is formulated.
As an application, we prove that discrete Lagrangian
dynamics on a Lie groupoid $G$ may be described in terms of Lagrangian implicit difference
equations of the corresponding cotangent groupoid $T^*G$. Other situations include finite
difference methods for time-dependent linear differential-algebraic equations and discrete
nonholonomic Lagrangian systems, as par\-ti\-cu\-lar examples.
\end{abstract}

\maketitle

\date{\today}

\begin{center}
\textit{Dedicated to Ernesto Lacomba on the occasion of his 65th birthday }
\end{center}
\section{Introduction}

Lie algebroids and groupoids have deserved  a lot of interest in recent years since these concepts generalize
the traditional framework  of tangent bundles and its discrete version, cartesian product of manifolds, to more
general situations. In particular,
it is well-known that the geometric description of the
Euler-Lagrange equations of a mechanical system determined by a
Lagrangian  function $L$, relies on the intrinsic geometry of the
tangent bundle $TQ$, the velocity phase space of a configuration
manifold $Q$. In the case when the Lagrangian $L$ is invariant
under the action of a Lie group $G$, the description, in this
case, relies on the geometry of the quotient space $TQ/G$ and the
equations describing the dynamics are called Lagrange-Poincar\'e
equations \cite{CMR}.  In this sense, Weinstein \cite{weinstein}
showed that the common geometric structure  of  the
Lagrange-Poincar\'e equations is essentially the same as the one
of the Euler-Lagrange equations, namely that of a Lie algebroid.
In the case of Lagrangian systems on the usual tangent bundles,
the tangent bundle carries a canonical Lie algebroid structure
which is given by the usual Lie algebra of vector fields on $Q$.
In the case of reduced lagrangian systems, we use the Atiyah
algebroid  $TQ/G\rightarrow Q/G$ to describe the evolution equations.

In \cite{MMM}, it is  described geometrically discrete Lagrangian and Hamiltonian Mechanics on Lie groupoids, in
particular, the type of equations analyzed  include  the classical discrete Euler-Lagrange equations, the
discrete Euler-Poincar\'e and discrete Lagrange-Poincar\'e equations. These results have applications  for the
construction of geometric integrators for continuous Lagrangian systems (reduced or not) \cite{Hair,Sanz}.

On the other hand, in \cite{Tul1,Tul2}, W.M. Tulczyjew proved that it is possible to interpret the ordinary
Lagrangian and Hamiltonian dynamics as Lagrangian submanifolds of suitable special symplectic manifolds. These
results were successfully extended to general Lie algebroids in \cite{LMM}.  More generally, in
\cite{MMT,MMT1,MMT2} implicit differential equations are considered. More precisely, an implicit differential
equation on a manifold $Q$ is a submanifold $D\subset TQ$. Given such a datum, the problem of integrability is
discussed, looking for a subset $S\subseteq D$ where, for any $v\in S$ there exists a curve $\gamma :I\to Q$
such that $\dot\gamma (0)=v$ and $\dot\gamma (t)\in D$ for any $t\in I$. An algorithm for extracting the
integrable part of an implicit differential equation is formulated. More precisely, if $D\subset TQ$ is an
implicit differential equation, the following sequence $\{ (C^k ,D^k )\}$ is considered
\[
\begin{array}{rl}
D^0=D, &  C^0=\tau _Q(D), \\
D^k= D^{k-1}\cap TC^{k-1}, & C^k=\tau _Q (D^k),
\end{array}
\]
where $\tau _Q:TQ\to Q$ is the canonical projection. Clearly, $D^k\subseteq D^{k-1}$ and $C^{k-1}\subseteq C^k$
for any $k$. Under some smoothness assumptions, it is proved in \cite{MMT1} that if the algorithm stabilizes,
that is, there exists $\overline{k}$ such that $D^{\bar k}=D^{{\bar k}-1}$, then $D^{\bar k}\subseteq D$ is the
(possibly empty) integrable part of $D$. The relation with Lagrangian mechanics is the following one. Given a
Lagrangian $L:TP\to \R$, one can consider the implicit differential equation $S_L\subseteq T(T^*P)$ using the
image of the differential of $L$, $dL(TP)\subset T^*(TP)$ and the canonical isomorphism $A_M:T^*(TP)\to T(T^*P)$. If we
apply the integrability algorithm, the first step will give us the solution of the Euler-Lagrange equations.
Indeed, if $(q^i,v^i)$ are local coordinates on $TP$ then
\[
\begin{array}{l}
\displaystyle S_L=\{ (q^i,\frac{\partial L}{\partial v^i};v^i,\frac{\partial L}{\partial q^i})\}, \\[5pt]
\displaystyle (S_L)^1 =\setdef{(q^i,\frac{\partial L}{\partial
v^i};v^i,\frac{\partial L}{\partial q^i})}{\dot{q}^i=v^i, \quad
\frac{d}{dt}\Big ( \frac{\partial L}{\partial v^i}\Big
)=\frac{\partial L}{\partial q^i}}.
\end{array}
\]

Taking as starting point the results of discrete mechanics on Lie groupoids
\cite{IMMM,MMM} and a remark on implicit differential equations on Lie algebroids
in \cite{IMMS}, in
this paper, submanifolds of a Lie groupoid $G$ are interpreted as systems of implicit
difference equations on $G$.
We study the problem of integrability of these implicit difference equations, proposing
an algorithm which extracts the integrable part (backward and/or forward) of this type of
systems and  analyze the geometric properties of them. In our
opinion, this approach offers greater conceptual clarity. Discrete Lagrangian systems
(unconstrained and nonholonomic) are studied within the framework of implicit difference
equations. Moreover, for the unconstrained case, the resulting implicit difference equation
is a Lagrangian submanifold, thus obtaining a discrete version of the results in
\cite{LMM}. Other examples include finite difference methods for time-dependent linear
differential algebraic equations.
The results about implicit difference equations are useful for systems more general than the
studied on this paper; in
particular, discrete optimal control theories, discrete systems with external constraints
... See for instance \cite{IMMM,MMS}.

\section{Preliminaries: Lie algebroids and Lie groupoids}

In this Section, we will recall the definition of a Lie groupoid and
some generalities about them are explained (for more details, see
\cite{Ma}).

A \emph{groupoid} over a set $M$ is a set $G$ together with the
following structural maps:
\begin{itemize}
\item A pair of surjective maps $\alpha: G \to M$, the \emph{source}, and $\beta: G \to M$, the \emph{target}.
These maps define the set of composable pairs
$$
G_{2}=\{(g,h) \in G \times G / \beta(g)=\alpha(h)\}.
$$
\item A \emph{multiplication} $m: G_{2} \to G$, to be denoted simply by $m(g,h)=gh$, such that
\begin{itemize}
\item[-] $\alpha(gh)=\alpha(g)$ and $\beta(gh)=\beta(h)$.
\item[-] $g(hk)=(gh)k$.
\end{itemize}
\item An \emph{identity section} $\epsilon: M \to G$ such that
\begin{itemize}
\item[-] $\epsilon(\alpha(g))g=g$ and $g\epsilon(\beta(g))=g$.
\end{itemize}
\item An \emph{inversion map} $i: G \to G$, to be denoted simply by $i(g)=g^{-1}$, such that
\begin{itemize}
\item[-] $g^{-1}g=\epsilon(\beta(g))$ and $gg^{-1}=\epsilon(\alpha(g))$.
\end{itemize}
\end{itemize}

A groupoid $G$ over a set $M$ will be denoted simply by the symbol
$G \rightrightarrows M$.

The groupoid $G \rightrightarrows M$ is said to be a \emph{Lie groupoid} if $G$ and $M$ are manifolds and all
the structural maps are differentiable with $\alpha$ and $\beta$ differentiable submersions. If $G
\rightrightarrows M$ is a Lie groupoid then $m$ is a submersion, $\epsilon$ is an injective immersion and $i$ is
a diffeomorphism. Moreover, if $x \in M$, $\alpha^{-1}(x)$ (resp., $\beta^{-1}(x)$) will be said the
\emph{$\alpha$-fiber} (resp., the \emph{$\beta$-fiber}) of $x$.

On the other hand, if $g \in G$ then the \emph{left-translation by
$g \in G$} and the \emph{right-translation by $g$} are the
diffeomorphisms
$$
\begin{array}{lll}
l_{g}: \alpha^{-1}(\beta(g)) \longrightarrow
\alpha^{-1}(\alpha(g))&; \; \;& h \longrightarrow
l_{g}(h) = gh, \\
r_{g}: \beta^{-1}(\alpha(g)) \longrightarrow \beta^{-1}(\beta(g))&;
\; \;& h \longrightarrow r_{g}(h) = hg.
\end{array}
$$
Note that $l_{g}^{-1} = l_{g^{-1}}$ and $r_{g}^{-1} = r_{g^{-1}}$.

A vector field $\tilde{X}$ on $G$ is said to be
\emph{left-invariant} (resp., \emph{right-invariant}) if it is
tangent to the fibers of $\alpha$ (resp., $\beta$) and
$\tilde{X}(gh) = (T_{h}l_{g})(\tilde{X}_{h})$ (resp., $\tilde{X}(gh)
= (T_{g}r_{h})(\tilde{X}(g))$, for $(g,h) \in G_{2}$.

Now, we will recall the definition of the \emph{Lie algebroid
associated with $G$}.

We consider the vector bundle $\tau: AG \to M$, whose fiber at a
point $x \in M$ is $A_{x}G = V_{\epsilon(x)}\alpha = Ker
(T_{\epsilon(x)}\alpha)$. It is easy to prove that there exists a
bijection between the space of sections of $\tau$, $\Gamma (\tau)$, and the set of
left-invariant (resp., right-invariant) vector fields on $G$. If $X$
is a section of $\tau: AG \to M$, the corresponding left-invariant
(resp., right-invariant) vector field on $G$ will be denoted
$\lvec{X}$ (resp., $\rvec{X}$), where
\begin{equation}\label{linv}
\lvec{X}(g) = (T_{\epsilon(\beta(g))}l_{g})(X(\beta(g))),
\end{equation}
\begin{equation}\label{rinv}
\rvec{X}(g) = -(T_{\epsilon(\alpha(g))}r_{g})((T_{\epsilon
(\alpha(g))}i)( X(\alpha(g)))),
\end{equation}
for $g \in G$. 

Using the above facts, we may introduce a Lie algebroid structure
$(\lcf\cdot , \cdot\rcf, \rho)$ on $AG$, which is defined by
\begin{equation}\label{LA}
\lvec{\lcf X, Y\rcf} = [\lvec{X}, \lvec{Y}], \makebox[.3cm]{}
\rho(X)(x) = (T_{\epsilon(x)}\beta)(X(x)),
\end{equation}
for $X, Y \in \Gamma(\tau)$ and $x \in M$. We recall that a
\emph{Lie algebroid} $A$ over a manifold $M$ is a real vector
bundle $\tau: A \to M$ together with a Lie bracket $\lcf\cdot,
\cdot\rcf$ on the space $\Gamma(\tau)$ of the global cross
sections of $\tau: A \to M$ and a bundle map, called \emph{the
anchor map}, such that if we also denote by $\rho:\Gamma(\tau)\to
{\mathfrak X}(M)$ the homomorphism of $C^\infty(M)$-modules
induced by the anchor map then
\[
\lcf X,fY\rcf=f\lcf X,Y\rcf + \rho(X)(f)Y,
\]
for $X,Y\in \Gamma(\tau)$ and $f\in C^\infty(M)$. The triple
$(A,\lcf\cdot,\cdot\rcf,\rho)$ is called a Lie algebroid over $M$
(see \cite{Ma}).

Next, we will present some examples of Lie groupoids and their
associated Lie algebroids which will be useful for our purposes.

Any Lie group $G$ is a Lie groupoid over $\{\mathfrak e \}$, the identity element of $G$, and the Lie algebroid
associated with $G$ is just the Lie algebra ${\mathfrak g}$ of $G$. On the other hand, given a manifold $M$, the
product manifold $M \times M$ is a Lie groupoid over $M$, called the \emph{pair or banal groupoid}, in the
following way: $\alpha$ is the projection onto the first factor and $\beta$ is the projection onto the second
factor, $m((x, y), (y, z)) = (x, z)$, for $(x, y), (y, z) \in (M \times M)_2$, $\epsilon(x) = (x, x)$, for all
$x \in M$, and $i(x, y) = (y, x)$. The Lie algebroid $A(M\times M)$ of this groupoid is isomorphic to the
tangent bundle $\tau_M:TM\to M$.

{\bf The cotangent groupoid.} (See \cite{CDW}, for details).
Let $G\rightrightarrows M$ be a Lie groupoid. If $A^\ast G$ is the dual bundle to
$AG$ then the cotangent bundle $T^\ast G$ is a Lie groupoid over
$A^\ast G$. The projections $\tilde\beta$ and $\tilde\alpha$, the
partial multiplication $\oplus _{T^\ast G}$, the identity section
$\tilde\epsilon$ and the inversion $\tilde{\i}$ are defined
as follows,
\begin{equation}\label{eq:cotangent:groupoid}
\kern-15pt
\begin{array}{l} \tilde\beta(\mu _g)(X)=\mu _g
((T_{\epsilon (\beta(g))} l_g )(X)), \mbox{ for }\mu _g\in
T^\ast _gG \mbox{ and }X\in A_{\beta(g)}G, \\[5pt] \tilde\alpha (\nu
_h)(Y)=\nu _h ((T_{\epsilon (\alpha (h))} r_h) (Y-
(T _{\epsilon (\alpha (h))}(\epsilon \circ \beta)) (Y))),\\[4pt]\kern160pt
\mbox{ for }\nu _h\in T^\ast _hG \mbox{ and }Y\in A_{\alpha(h)}G ,
\\[5pt] (\mu _g\oplus _{T^\ast G}\nu _h)(T_{(g,h)}m (X_g,Y_h))=\mu
_g(X_g)+\nu _h (Y_h),\\[4pt]\kern160pt\mbox{ for }(X_g,Y_h)\in
T_{(g,h)}G_{2},
\\[5pt] \tilde\epsilon(\mu _x)(X_{\epsilon (x)})=\mu _x(X_{\epsilon
(x)}- (T_{\epsilon (x)}(\epsilon \circ \alpha )) (X_{\epsilon
(x)}))),\\[4pt] \kern160pt \mbox{ for }\mu _x\in A^\ast _xG \mbox{
and }X_{\epsilon (x)}\in T_{\epsilon (x)}G ,\\[5pt] \tilde{\i}
(\mu _g)(X_{g^{-1}})=-\mu _g((T_{g^{-1}} i )(X_{g^{-1}})),
\mbox{ for }\mu _g\in T^\ast _gG\mbox{ and }X_{g^{-1} }\in
T_{g^{-1}}G.
\end{array}
\end{equation}
Note that $\tilde\epsilon( A^\ast G )$ is just the conormal bundle of $M\cong \epsilon (M)$ as a submanifold of
$G$. In addition, $A^*G$ is endowed with a linear Poisson structure. For this Poisson structure on $A^*G$,
$\tilde\alpha:T^*G \to A^*G$ is a Poisson map and $\tilde\beta:T^*G \to A^*G$ is an anti-Poisson (for more details,
see \cite{CDW,Ma}).

\section{Discrete dynamics in implicit form}

Motivated by the studies of implicit differential equations in
\cite{MMT,MMT1,MMT2} and its generalization to Lie algebroids in
\cite{IMMS}, we introduce the notion of implicit difference
equations on Lie groupoids and the corresponding notion of
integrability.

\begin{definition}
An \emph{implicit difference equation} on a Lie groupoid $G
\rightrightarrows M$ is a submanifold $E$ of $G$.

An \textbf{admissible sequence} on the Lie groupoid $G\gpd M$ is a mapping $\gamma_G: I\cap \Z\longrightarrow G$
such that $\beta(\gamma_G(i))=\alpha(\gamma_G(i+1))$ for all $i,i+1\in I\cap \Z$. Here $I$ is an interval on
$\R$.

A \emph{solution} of an implicit difference equation $E\subset G$ is an admissible sequence $\gamma_G: I\cap
\Z\longrightarrow G$  on $G$ such that $\gamma _G(i)\in E$, for all $i\in I\cap\Z$.
\end{definition}

Let $E$ be a submanifold of the Lie groupoid $G\gpd M$. Then $E$ is said to be
\begin{enumerate}
\item[1)] \emph{forward integrable at} $g\in E$  if there is a solution
$\gamma_G: \Z^+ \longrightarrow E\subseteq G$  with $\gamma_G(0)=g$.
Here, $\Z^+=\{n\in \Z\; |\; n\geq 0\}$.
\item[2)] \emph{backward integrable at} $g\in E$ if there is a solution
$\gamma_G: \Z^- \longrightarrow E\subseteq G$ with $\gamma_G(0)=g$.
Here, $\Z^-=\{n\in \Z\; |\; n\leq 0\}$.
\item[3)] \emph{integrable at} $g\in E$ if if there is a solution
$\gamma_G: \Z \longrightarrow E\subseteq G$ with $\gamma_G(0)=g$, that is,
it is backward and forward integrable at $g\in E$.
\end{enumerate}
If these conditions hold for all $g$, we say that $E$ is forward integrable, backward integrable or
integrable, respectively.

\begin{proposition}\label{prop:first:charact}
Let $E$ be a submanifold of the Lie groupoid $G\gpd M$. Then,
\begin{enumerate}
\item[1)] $E$ is forward integrable if and only if for each $g\in E$ exists at least
an $h\in E$ such that $(g, h)\in G_2$ or, equivalently, $E\subseteq
\beta^{-1}(\alpha(E))$.
\item[2)] $E$ backward integrable if and only if for each $g\in E$ exists at
least an $h'\in E$ such that $(h', g)\in G_2$ or, equivalently, $E\subseteq
\alpha^{-1}(\beta(E))$.
\item[3)] $E$ is integrable if and only if $E\subseteq
\alpha^{-1}(\beta(E))\cap\beta^{-1}(\alpha(E))$.
\end{enumerate}
\end{proposition}
\begin{proof}
1) Suppose that $E$ is forward integrable. Then, for all $g\in E$ there is an admissible
sequence $\gamma _G:\Z^+\to E$ such that $\gamma _G(0)=g$. So, we have that $h=\gamma _G(1)\in E$ satisfies
$(g,h)\in G_2$.

Conversely, if $g\in E$ there exists a $g_1\in E$ such that $(g,g_1)\in G_2$. Since $g_1\in E$ and using again
the hypothesis, there is $g_2\in E$ satisfying $(g_1,g_2)\in G_2$. Therefore, we can construct a sequence
$\{g_i\}_{i\in\Z^+}\subseteq E$ where $g_0=g$ and $(g_i,g_{i+1})\in G_2$. As a consequence, the sequence $\gamma
_G:\Z^+\to E$, $i\mapsto \gamma_G(i)=g_i$, is a solution of $E$ with $\gamma_G(0)=g$.

With a similar argument we deduce 2), and 3) is a direct consequence of 1) and 2).
\end{proof}
\begin{example}{\rm 
Let $M$ be a manifold and $\varphi \in C^\infty (M)$ be a diffeomorphism on $M$.
Consider the implicit difference equation $E_\varphi$ of the pair groupoid $M\times M$
given by
\[
E_\varphi =\{ (x,\varphi (x))\in M\times M\, | \, x\in M\}.
\]
Then, if $(x_0,\varphi (x_0))\in E_\varphi$, the map $\gamma :\Z \to E_\varphi$ defined by
\[
\gamma (i)=(\varphi ^i (x_0),\varphi ^{i+1} (x_0)), \quad \mbox{ for }i\in \Z,
\]
where $\varphi ^{-i}=(\varphi^{-1})^{i}$, is a solution of the implicit difference equation $E_\varphi$. Therefore, $E_\varphi$ is integrable.

This example can be generalized to bisections on Lie groupoids. We recall that
a bisection on a Lie groupoid is a map $\sigma :M\to G$ which is a section of
$\alpha$ and such that $\beta\circ \sigma$ is a diffeomorphism. Given a bisection,
$\sigma$, if $E_\sigma =\sigma (M)$ is the associated implicit difference equation on $G$  then it is integrable, since for any $g=\sigma(x)\in E_\sigma$, the sequence
\[
\gamma (i) = \sigma ( (\beta \circ \sigma )^{i} (x) ),\quad \mbox{ for }i\in \Z,
\]
is a solution with $\gamma (0)=g$.}
\end{example}

The sets
\begin{eqnarray*}
\tilde{E}_f&=&\{g\in E\;|\; E \hbox{ is {forward integrable at} }
g\}, \\
\tilde{E}_b&=&\{g\in E\;|\;  E \hbox{ is {backward integrable at}
}
g\},\\
\tilde{E}_{fb}&=&\{g\in E\;|\; E \hbox{ is {integrable at} } g\},
\end{eqnarray*}
are called the \emph{forward integrable, backward integrable} and \emph{integrable parts} of $E$, respectively.
The implicit difference equation $E$ is integrable (resp. forward integrable or backward integrable) if
$E=\tilde{E}_{fb}$ (resp. $E=\tilde{E}_f$ or $E=\tilde{E}_b$). Note that if $E\subseteq E'$ then
$\tilde{E}_{f}\subset \tilde{E}'_{f}$, $\tilde{E}_{b}\subset \tilde{E}'_{b}$
 and $\tilde{E}_{fb}\subset \tilde{E}'_{fb}$.

\begin{proposition}
Let $E$ be a submanifold of a Lie groupoid $G\gpd M$. Assume that $\tilde{E}_{fb}$, the integrable part of $E$,
is a submanifold. Then, $\tilde{E}_{fb}$ is integrable.
\end{proposition}
\begin{proof}
If $g\in \tilde{E}_{fb}$ then there exists an admissible sequence $\gamma_G: \Z
\to E\subseteq G$ with $\gamma_G(0)=g$. We have that $\gamma_G(\Z)\subseteq \tilde{E}_{fb}$. Indeed,
for all $\tilde{n}\in \Z$ there is an admissible sequence $\gamma_{G, \tilde{n}}: \Z\to E$ given by
\[
\gamma_{G, \tilde{n}}(i)= \gamma_G(i+\tilde{n}),\quad \mbox{ for all }i\in\Z.
\]
\end{proof}
A similar result holds for $\tilde{E}_f$ and $\tilde{E}_b$. Thus,
$\tilde{E}_f$ and $\tilde{E}_b$ are forward and
backward integrable, respectively.

Now, it is easy to prove the following
\begin{proposition}\label{aa2}
Let $E$ be a submanifold of a Lie groupoid $G\gpd M$.
\begin{enumerate}
\item[(i)] If $\tilde{E}_f$ is the forward integrable part of $E$ and $E'$ is a submanifold of $G$ such that
$\tilde{E}_f\subseteq E'\subseteq E$, then $\tilde{E}_f$ is the forward integrable part of $E'$. \item[(ii)] If
$\tilde{E}_b$ is the backward integrable part of $E$ and $E'$ is a submanifold of $G$ such that
$\tilde{E}_b\subseteq E'\subseteq E$, then $\tilde{E}_b$ is the backward integrable part of $E'$. \item[(iii)]
If $\tilde{E}_{fb}$ is the  integrable part of $E$ and $E'$ is a submanifold of $G$ such that
$\tilde{E}_{fb}\subseteq E'\subseteq E$, then $\tilde{E}_{fb}$ is the integrable part of $E'$.
\end{enumerate}
\end{proposition}


From Proposition \ref{prop:first:charact}, we have the following obvious
relations
\begin{eqnarray*}
 \tilde{E}_f&\subseteq&
\beta^{-1}(\alpha(\tilde{E}_f))\subseteq \beta^{-1}(\alpha(E))\\
\tilde{E}_b&\subseteq&
\alpha^{-1}(\beta(\tilde{E}_b))\subseteq \alpha^{-1}(\beta(E))\\
\tilde{E}_{fb}&\subseteq&
\beta^{-1}(\alpha(\tilde{E}_{fb}))\cap\alpha^{-1}(\beta(\tilde{E}_{fb}))\subseteq
\beta^{-1}(\alpha({E}))\cap\alpha^{-1}(\beta({E}))
\end{eqnarray*}
Therefore,
\begin{eqnarray*}
 \tilde{E}_f&\subseteq&
E\cap \beta^{-1}(\alpha(E))\\
\tilde{E}_b&\subseteq&
E\cap \alpha^{-1}(\beta(E))\\
\tilde{E}_{fb}&\subseteq& E\cap
\beta^{-1}(\alpha({E}))\cap\alpha^{-1}(\beta({E}))
\end{eqnarray*}

Propositions \ref{prop:first:charact} and \ref{aa2} suggest a method for extracting the integrable part of an
implicit difference equation under some regularity hypothesis.

Let $E\subset G$ be an implicit difference equation of a Lie groupoid
$G\rightrightarrows M$.

{\bf For forward integrable part}. Consider the sequence of sets of
$M$
\[
C_f^0=\alpha(E), \quad C_f^1=\alpha(E\cap \beta^{-1}(C_f^0)),
\ldots, C_f^k=\alpha(E\cap \beta^{-1}(C_f^{k-1})), \ldots
\]
and the sequence subsets of $G$
\[
E_f^0=E, \quad E_f^1=E\cap \beta^{-1}(C_f^0), \ldots, E_f^k=E\cap
\beta^{-1}(C_f^{k-1}), \ldots
\]
Observe that $\alpha(E_f^k)=C_f^k$,  $C_f^k\subseteq C_f^{k-1}$ and
$ E_f^k\subseteq E_f^{k-1}$.

{\bf For backward integrable part}.  Consider the sequence of sets of
$M$
\[
D_b^0=\beta (E), \quad D_b^1=\beta(E\cap \alpha^{-1}(D_b^0)),
\ldots, D_b^k=\beta(E\cap \alpha^{-1}(D_b^{k-1})), \ldots
\]
and the sequence of subsets of $G$
\[
E_b^0=E, \quad E_b^1=E\cap \alpha^{-1}(D_b^0), \ldots, E_b^k=E\cap
\alpha^{-1}(D_b^{k-1}), \ldots
\]
Also then $ D_b^k \subseteq D_b^{k-1}$ and $E_b^k \subseteq E_b^{k-1}$.

{\bf For integrable part}.  Consider the sequence of sets of $M$
\begin{eqnarray*}
C_{fb}^0&=&\alpha(E),\quad  D_{fb}^0=\beta(E),\\
C_{fb}^1&=&\alpha(E\cap \alpha^{-1}(D_{fb}^0)\cap
\beta^{-1}(C_{fb}^0)), \quad D_{fb}^1=\beta(E\cap
\alpha^{-1}(D_{fb}^0)\cap \beta^{-1}(C_{fb}^0))
\\
 C_{fb}^k&=&\alpha(E\cap \alpha^{-1}(D_{fb}^{k-1})\cap
\beta^{-1}(C_{fb}^{k-1})), \quad D_{fb}^{k}=\beta(E\cap \alpha^{-1}(D_{fb}^{k-1})\cap \beta^{-1}(C_{fb}^{k-1}))
\ldots
\end{eqnarray*}
and the sequence of subsets of $G$
\[
\begin{array}{l}
E_{fb}^0=E,  E_{fb}^1=E\cap \alpha^{-1}(D_{fb}^0)\cap \beta^{-1}(C_{fb}^0), \ldots\\[7pt]
  E_{fb}^k=E\cap \alpha^{-1}(D_{fb}^{k-1}))\cap \beta^{-1}(C_{fb}^{k-1}), \ldots
\end{array}
\]
Then $C_{fb}^k \subseteq C_{fb}^{k-1}$,
$D_{fb}^k \subseteq D_{fb}^{k-1}$ and $E_{fb}^k\subseteq E_{fb}^{k-1}$.



It may happen that after a finite number of steps, the sets in two consecutive steps in the different sequences
are equal:
\begin{eqnarray*}
&C^{k_f}_f=C^{{k_f}+1}_f=\bar{C}_f,\quad D^{k_b}_b=D^{k_b+1}_b=\bar{D}_b,&\\
&C^{k_{fb}}_{fb}=C^{k_{fb}+1}_{fb}=\bar{C}_{fb},\quad
D^{k_{fb}}_{fb}=D^{k_{fb}+1}_{fb}=\bar{D}_{fb}.&
\end{eqnarray*}
Then
\begin{eqnarray*}
&E^{k_f+1}_f=E^{k_f+2}_f=\bar{E}_f,\quad E^{k_b+1}_b=E^{k_b+2}_b=\bar{E}_b,&\\
&E^{k_{fb}+1}_{fb}=E^{k_{fb}+2}_{fb}=\bar{E}_{fb}.&
\end{eqnarray*}
It is straightforward to see that if this happens then the subsequent steps in the different
sequences are all equal.

Moreover, observe that $\alpha(\bar{E}_f)=\bar{C}_f$, $\beta(\bar{E}_b)=\bar{D}_b$,
$\alpha(\bar{E}_{fb})=\bar{C}_{fb}$ and $\beta(\bar{E}_{fb})=\bar{D}_{fb}$.

\begin{theorem}
The forward integrable part of $E$ is $\bar{E}_f$, the backward integrable
part is $\bar{E}_b$ and the integrable part of
$E$ is $\bar{E}_{fb}$.
\end{theorem}
\begin{proof}
Let $g\in \bar{E}_f=E^{k_f}_f=E^{k_f+1}_f$. Then, since $E^{k_f+1}_f= E\cap \beta^{-1}(C^{k_f}_f)$, we have that
$\beta (g)\in C^{k_f}_f=\alpha(\bar{E}_f)$, that is, there exists $h\in E^{k_f+1}_f= E^{k_f}_f= \bar{E}_f$ such
that $\beta (g)=\alpha (h)$. Therefore, using Proposition \ref{prop:first:charact}, we have that $E$ is forward
integrable at $g$.

Conversely, suppose that $E$ is forward integrable at $g$. Then, there exists an admissible sequence $\gamma _G
:\Z^+\to E\subseteq G$ such that $\gamma _G(0)=g$. Denote $\gamma _G(i)$ by $g_i$, for all $i\in\Z^+$. Then,
$\beta (g_{k_f})=\alpha (g_{{k_f}+1})$ and $g_{k_f}\in E\cap \beta ^{-1}(C^0_f)=E^1_f$, which implies that
$\alpha (g_{k_f})\in C^1_f$. Now, since $\beta (g_{{k_f}-1})=\alpha (g_{k_f})$, we deduce that $g_{{k_f}-1}\in
E^2_f$. If we repeat this argument ${k_f}$ times, we conclude that $g=g_0\in E_f^{k_f+1}=E^{k_f+2}_f=\bar{E}_f$.

The corresponding results for backward integrability and integrability are proved in a similar way.
\end{proof}

\begin{remark}\label{rem:generalization}
We have to stress that the concept of implicit difference equation
can be to generalized to the setting of a pair of manifolds $G$ and
$M$ endowed with two surjective submersions $\alpha ,\beta:G\to M$.
The corresponding notions of admissible sequences, solutions,
integrable equations and the integrability algorithm can be
straightforwardly generalized. This more general approach will be
useful in Section \ref{sub:NH}.
\end{remark}

\begin{remark}
For some concrete applications (see, for instance, the next example) instead of a unique subset $E$ of the Lie groupoid $G\rightrightarrows M$,  we have a sequence of subsets $\{E_k\}_{0\leq k\leq \infty}$. A generalization of this type is also useful for discretization of time-dependent systems.  It is not hard to extend our setting for this situation and also to extract the integrable part of a sequence of this type. For example, the forward integrable part  is extracted using the following procedure. Consider the $k$-sequence of sets of $M$ defined for $k\geq 0$:
\[
[C_k]^0_f=\alpha(E_k), \quad [C_k]^1_f=\alpha(E_{k}\cap \beta^{-1}([C_{k+1}]^0_f)), \ldots, [C_k]^i_f=\alpha(E_{k}\cap \beta^{-1}([C_{k+1}]^{i-1}_f)),
\]
and the subsets
\[
[E_{k}]^0_f=E_{k}, \quad [E_{k}]^1_f=E_{k}\cap \beta^{-1}([C_{k+1}]^0_f), \ldots, [E_{k}]^i_f=E_{k}\cap \beta^{-1}([C_k]^{i-1}_f).
\]
If the algorithm stabilizes in a sequence of subsets
$\{\overline{[E_k]}_f\}$ we will call this set the backward integrable part of the sequence $\{E_k\}$.
\end{remark}

\begin{example}{\bf Finite difference methods for time-dependent linear diffe\-rential-algebraic equations}
{\rm 
The proposed algorithm is also useful for numerical methods for differential-algebraic equations (DAEs). These systems of equations appear in a great number of applications: mechanical systems with constraints, electrical networks, chemical reacting systems... See, for instance,  \cite{CeEt,PeLo}.  These types of systems take the form $F(t, x, \dot{x})=0$ including as a particular example the explicit differential equations.
The system that we will study are the so-called time-dependent linear differential-algebraic equations:
\[
A(t)\dot{x}+B(t)x=b(t), \qquad A(t), B(t)\in {\mathcal L}(\R^n, \R^n),\;  b(t)\in \R^n,\;  t\in \R
\]

Consider an explicit Euler approximation to this system
\[
A_k\frac{x_{k+1}-x_k}{h}+B_kx_k=b_k
\]
where $h>0$ is the step size and $A_k=A(t_0+kh)$, $B_k=B(t_0+kh)$ and $b_k=b(t_0+kh)$ and $x_0, t_0$ initial conditions.
Define the subsets of the pair groupoid $\R^n\times \R^n$
\[
E_k=\{(x_k, x_{k+1})\in \R^n\times \R^n \; \mid \; A_k\frac{x_{k+1}-x_k}{h}+B_kx_k=b_k\}
\]
We consider $V_k$ the image of the linear map associated with $A_k$. Then, $V_k=\R^n$
(in this case, $x_{k+1}=x_k+hA_k^{-1}(b_k-B_kx_k)$) or there exists a matrix $Q_k\in
{\mathcal L}(\R^n, \R^n)$ such that the kernel of its linear map is just $V_k$, i.e.,
$Q_kA_k={\mathbf 0}$. Then,
one may easily prove that
\[
[C_k]^0_f=\{x\in \R^n\; \mid\; Q_{k}B_{k}x-Q_{k}b_{k}=0\}\; .
\]
Therefore,
\begin{eqnarray*}
[E_{k}]^1_f&=&\{(x_k, x_{k+1})\in \R^n\times \R^n \; \mid \; A_k\frac{x_{k+1}-x_k}{h}+B_kx_k=b_k,\\
&& Q_{k+1}B_{k+1}x_{k+1}-Q_{k+1}b_{k+1}=0\}\\
&=&\{(x_k, x_{k+1})\in \R^n\times \R^n \; \mid \; Q_{k+1}B_{k+1}x_{k+1}-Q_{k+1}b_{k+1}=0,\ \\
&&(A_k+Q_{k+1}B_{k+1})x_{k+1}=(A_k-hB_k)x_k+hb_k+Q_{k+1}b_{k+1},\}\; .
\end{eqnarray*}
It is possible to continue the algorithm, considering $\alpha([E_{k}]^1_f)=[C_k]^1_f$, but in many cases of interest (as, for instance, discretizations of index one DAEs) we have that the matrix $A_k+Q_{k+1}B_{k+1}$ is regular, then
\[
x_{k+1}=(A_k+Q_{k+1}B_{k+1})^{-1}\left((A_k-hB_k)x_k+hb_k+Q_{k+1}b_{k+1}\right)\; .
\]
In this case, $\alpha([E_{k}]^1_f)=[C_k]^1_f=[C_k]^0_f$ and the algorithm stops in this step (see \cite{RaRh} for more details).}
\end{example}
\section{Implicit systems defined by discrete Lagrangian functions}

\subsection{Implicit systems defined by discrete Lagrangian functions. Lagran\-gian
submanifolds}

Let $L:G\to \R$ be a discrete Lagrangian on a Lie groupoid $G\gpd M$. Then, $S_L=dL(G)\subseteq T^\ast G$ is
clearly a Lagrangian submanifold of the symplectic manifold $T^\ast G$ and $S_L$ is also an implicit difference
equation on the cotangent groupoid $T^*G\gpd A^*G$. If we apply the forward integrability algorithm to $S_L$,
then there is an interesting relation with the discrete Euler-Lagrange equations for $L$ (see \cite{MMM,Stern}
and the Appendix \ref{A} for more details on discrete Lagrangian Mechanics on Lie groupoids), which we describe in this
section.
\begin{proposition}\label{forward-DEL}
Let $G\gpd M$ be a Lie groupoid and $L:G\to \R$ be a discrete Lagrangian on $G$. The first step, $(S_L)^1_f$, to
obtain the forward integrable part of $S_L=dL(G)$ are the points $dL(g)\in S_L$ such that there exists an
element $h\in G$, satisfying $(g, h)\in G_2$ and, in addition, $(g,h)$ is a solution of the discrete
Euler-Lagrange equations, that is,
\[
d\left[L\circ l_{g}+L\circ r_{h}\circ i\right](\epsilon (x))_{|
A_x G}=0,
\]
where $x=\beta (g)=\alpha (h)$.

As a consequence, if $\tilde\beta$ (resp. $\tilde\alpha$) is the target (resp. source) map of the cotangent
groupoid $T^*G\gpd A^*G$ defined in \eqref{eq:cotangent:groupoid}, then
\[
\begin{array}{rcl}
(S_L)^1_f&=&\{ (g,h)\in G_2\, | \, (g,h) \mbox{ is a solution of the discrete Euler-Lagrange eqns.}\}\\
&=&\{ (g,h)\in G_2\, | \, \tilde\beta (dL(g))=\tilde\alpha (dL(h))\}
\end{array}
\]
\end{proposition}
\begin{proof}
Using the definition of the source map $\tilde\alpha:T^*G\to A^*G$ (see \eqref{eq:cotangent:groupoid}) and the
fact that
\[
X-T_{\epsilon (x)}(\epsilon \circ \beta)(X)= - T_{\epsilon (x)}i
(X),\mbox{ for } X\in A_xG,
\]
we have that
\[
C_f^0=\tilde\alpha (S_L)=\{a^*\in A^*G\, |\, a^*=-d(L\circ
r_h\circ i)(\epsilon(\alpha (h)))\in A^*_{\alpha(h)}G, \hbox{ for
some } h\in G\}.
\]
Therefore, $(S_L)^1_f=S_L\cap (\tilde\beta)^{-1}(C_f^0)$ is
the set of points $dL(g)\in S_L$ such that there exists $h\in G$,
with $(g, h)\in G_2$ and
\[
d\left[L\circ l_{g}+L\circ r_{h}\circ i\right](\epsilon (x))_{|
A_x G}=0.
\]
\end{proof}

Let $L:G\to \R$ be a hyperregular discrete Lagrangian function (see Appendix \ref{A}). This implies that the discrete
Lagrangian evolution operator $\Upsilon _L:G\to G$ is given by $\Upsilon _L= (\F ^-L)^{-1}\circ (\F^+L)$. Then,
if $g\in G$, the map $\Upsilon _{T^*G}:\Z \to S_L=dL(G)\subseteq T^*G$ defined by
\[
\Upsilon _{T^*G}(i)=dL (\Upsilon _L^i (g)), \quad \mbox{ for }i\in \Z,
\]
is a solution of the implicit difference equation $S_L$ and $\Upsilon _{T^*G}(0)=dL (g)$.

Thus, we have proved the following result.
\begin{corollary}\label{cor:regularity}
Let $G\gpd M$ be a Lie groupoid and $L:G\to \R$ be a discrete hyperregular Lagrangian function on $G$. Then, the
implicit difference equation $S_L=dL(G)$ is integrable.
\end{corollary}
Under the same hypotheses as in Corollary \ref{cor:regularity}, we may define the map
\[
\begin{array}{l}
\Phi :G \stackrel{(Id,\Upsilon _L)}{\longrightarrow}G\times G
\stackrel{(\F ^-L,\F ^-L)}{\longrightarrow}A^\ast G\times A^\ast G\\
\quad\quad g\longmapsto (g,\Upsilon _L(g))\longmapsto (\F^-L(g),
\F^+L(g))
\end{array}
\]
A direct computation, using \eqref{eq:Leg:transform}, shows that $\Phi (G)=(\tilde\alpha,\tilde\beta)(S_L)$.

On the other hand, if $\tilde{\Upsilon}_L=\F^+L\circ (\F^-L)^{-1}
:A^\ast G\to A^\ast G$ is the discrete Hamiltonian evolution
operator then one observes that
\[
graph\, \tilde{\Upsilon}_L = \Phi (G)=(\tilde\alpha,\tilde\beta)(S_L).
\]

Using that $S_L$ is a Lagrangian submanifold of $T^*G$ and
that $(\tilde\alpha,\tilde\beta):T^*G\to A^*G\times
\overline{A^*G}$ is a Poisson map, $\overline{A^*G}$ denoting $A^*G$ endowed with the
linear Poisson structure changed of sign, then $(\tilde\alpha,\tilde\beta
)(S_L)$ is a coisotropic submanifold of $A^*G\times \overline{A^*G}$ (this is
a particular case of \cite[Corollary 2.2.5]{weinstein:coisotropic}).

From the previous discussion, we can conclude the following result (see also \cite{weinstein}).
\begin{proposition}
Let $G\gpd M$ be a Lie groupoid and $L:G\to \R$ be a discrete hyperregular Lagrangian on $G$. Then, the discrete
Hamiltonian evolution operation $\tilde{\Upsilon}_L$ preserves the Poisson structure on $A^*G$.
\end{proposition}

\subsection{Example: the pair groupoid} Let $Q$ be a manifold and
$L:Q\times Q\to \R$ be a discrete Lagrangian. We know that $dL(Q\times Q)$ is a Lagrangian submanifold of
$T^\ast (Q\times Q)$. On the other hand, $T^\ast (Q\times Q)\gpd T^\ast Q$ is a symplectic groupoid with the
canonical symplectic structure of $T^\ast (Q\times Q )$. The structural maps of $T^*(Q\times Q)\gpd T^*Q$ are
\[
\begin{array}{l}
\tilde{\beta}:T^\ast (Q\times Q)\to T^\ast Q,\quad (\gamma
_{q_0},\gamma _{q_1})\mapsto \gamma _{q_1},\\
\tilde{\alpha}:T^\ast (Q\times Q)\to T^\ast Q,\quad (\gamma
_{q_0},\gamma _{q_1})\mapsto -\gamma _{q_0},\\
\tilde{m}:T^\ast (Q\times Q)_2\to T^\ast Q,\quad \Big ( (\gamma
_{q_0},\gamma _{q_1}),(-\gamma _{q_1},\gamma _{q_2})\Big )\mapsto
(\gamma _{q_0},\gamma _{q_2}),\\
\tilde{i}:T^\ast (Q\times Q)\to T^\ast (Q\times Q),\quad (\gamma
_{q_0},\gamma _{q_1})\mapsto (-\gamma _{q_1},-\gamma _{q_0}),\\
\tilde{\epsilon}:T^\ast Q\to T^\ast (Q\times Q),\quad \gamma
_q\mapsto (\gamma _q,-\gamma _q).
\end{array}
\]
Moreover,
\[
\begin{array}{l}
\tilde{\beta}(dL(q_0,q_1))=D_2L(q_0,q_1),\\
\tilde{\alpha}(dL(q_1,q_2))=-D_1L(q_1,q_2),
\end{array}
\]
for $(q_0,q_1),(q_1,q_2)\in Q\times Q$. Thus, the discrete Euler-Lagrange
equations for the points $(q_0,q_1),(q_1,q_2)$, that is
$D_2L(q_0,q_1)+D_1L(q_1,q_2)=0$, are equivalent to
\[
\tilde{\beta}(dL(q_0,q_1)) =\tilde{\alpha}(dL(q_1,q_2)).
\]
Therefore, we obtain the following result (which can also be seen as a consequence of Proposition \ref{forward-DEL})
\begin{corollary}
Let $Q$ be a smooth manifold and $L:Q\times Q\to \R$ be a discrete Lagrangian function. Then, there exists a
bijective correspondence between composable points in $dL(Q\times Q)$ and solutions of the discrete
Euler-Lagrange equations.
\end{corollary}


\begin{remark}
{\rm $T^\ast (Q\times Q) \gpd T^\ast Q$ is a Lie groupoid which is
symplectomorphic to the pair groupoid $T^\ast Q\times T^ *Q\gpd
T^\ast Q$ when we consider the symplectic form $(-\Omega _Q,\Omega
_Q)$ on $T^\ast Q\times T^\ast Q$.

The isomorphism is given by
\begin{equation}\label{Psi}
\begin{array}{rcl}
\Psi:T^\ast (Q\times Q)&\to &T^\ast Q\times T^\ast Q\\
(\gamma _{q_0},\gamma _{q_1})&\mapsto &(-\gamma _{q_0},\gamma
_{q_1})
\end{array}
\end{equation}
Note that $(\Psi \circ
dL)(q_0,q_1)=(-D_1L(q_0,q_1),D_2L(q_0,q_1))$. }
\end{remark}
Now, suppose that $L:Q\times Q\to \R$ is hyperregular. Then 
\begin{equation}\label{operador}
\Upsilon _L =(\F^-L)^{-1}\circ \F ^+L:Q\times Q\to
Q\times Q
\end{equation}
is a discrete Lagrangian evolution operator. In addition, we have that the image of the map
\[
\begin{array}{l}
Q\times Q \stackrel{(Id,\Upsilon _L)}{\to}(Q\times Q)\times (Q\times
Q)\stackrel{\F^-L\times \F^-L}{\to}T^\ast Q\times T^\ast Q\\
(q_0,q_1)\mapsto ((q_0,q_1),\Upsilon _L (q_0,q_1))\mapsto
(-D_1L(q_0,q_1),D_2L(q_0,q_1))
\end{array}
\]
is just the Lagrangian submanifold $\Psi \circ dL(Q\times Q)$.

On the other hand, we can consider the discrete Hamiltonian
evolution operator $\tilde{\Upsilon}_L =\F ^+L\circ (\F^-L)^{-1}:T^\ast
Q\to T^\ast Q$. A direct computation proves that
\[
\Phi(Q\times Q)=(\Psi \circ dL)(Q\times Q)=(Id, \tilde{\Upsilon}_L)(T^\ast Q).
\]
We remark that, since $\tilde{\Upsilon}_L:T^\ast Q\to T^\ast Q$ is a
symplectomorphism, we deduce $(Id, \tilde{\Upsilon}_L)(T^\ast Q)$ is a
Lagrangian submanifold of $(T^\ast Q\times T^\ast Q,(-\Omega
_Q,\Omega _Q))$.

\subsection{Example: a discrete singular Lagrangian}
Consider the discretization of the degenerate Lagrangian $L:
T\R^2\longrightarrow \R$  given by $L(x, y, \dot{x},
\dot{y})=\frac{1}{2}\dot{x}^2+\frac{1}{2}x^2y$ (see \cite{Ca}).
Take, for instance,  a consistent discretization of this
lagrangian:
\[
L^h_d: \R^2\times \R^2\longrightarrow \R, \quad L^h_d(x_1, y_1, x_2,
y_2)=\frac{1}{2}\left(\frac{x_2-x_1}{h}\right)^2+\frac{1}{2}x_1^2y_1
\]
Then $S_{L_d^h}\subset T^*(\R ^2\times \R^2)$ defined by
\[
S_{L_d^h}=\left\{(x_1, y_1, x_2, y_2; -\frac{x_2-x_1}{h^2}+x_1y_1,
\frac{1}{2}x_1^2, \frac{x_2-x_1}{h^2},0)\; |\; (x_1,
y_1, x_2, y_2)\in \R^4\right\}
\]
The algorithm produces the sequence:
\begin{eqnarray*}
C^1_f&=&\tilde{\alpha}(S_{L_d^h})=\{ (x, y;
\frac{\tilde{x}-x}{h^2}-xy, -\frac{1}{2}x^2)\; | (x, y,
\tilde{x})\in \R^3)\} \subset  T^*\R^2\\
\left(S_{L_d ^h}\right)^1_f&=&\left\{(x_1, y_1, 0, y_2;  \frac{x_1}{h^2}+x_1y_1, \frac{1}{2}x_1^2,
 -\frac{x_1}{h^2},0)\; |\; (x_1, y_1, y_2)\in \R^3\right\} \\
C^2_f&=& \tilde{\alpha}(\left(S_{L_d^h}\right)^1_f)=\{ (x, y;
-\frac{x}{h^2}-xy, -\frac{1}{2}x^2)\; |\, (x, y)\in \R^2)\} \subset  T^*\R^2\\
\left(S_{L_d^h}\right)^2_f&=&\left\{(0, y_1, 0, y_2; 0, 0, 0, 0)\;
|\; (y_1, y_2)\in \R^2\right\}\\ C^3_f&=&
\tilde\alpha(\left(S_{L_d^h}\right)^2_f)=\{ (0, y; 0, 0)\; |
y\in \R)\} \subset  T^*\R^2
\end{eqnarray*}

Since $\left(S_{L_d^h}\right)^3_f=\left(S_{L_d^h}\right)^2_f$ then
the algorithm terminates with
$\left(\overline{S_{L_d^h}}\right)_f=\left(S_{L_d^h}\right)^2_f$, which is the
forward integrable part of the implicit difference equation
$S_{L_d^h}$.

\subsection{Implicit systems defined by discrete nonholonomic Lagrangian systems}
\label{sub:NH}

Let $(L,{\mathcal M}_c,{\mathcal D}_c)$ be a discrete nonholonomic
Lagrangian system. Consider the submanifold $S_{(L,{\mathcal
M}_c)}\subseteq T^*G$ given by $S_{(L,{\mathcal M}_c)}=dL({\mathcal
M}_c)$, that is,
\[
S_{(L,{\mathcal M}_c)}=\set{dL(i_{{\mathcal M}_c}(g))}{g\in{\mathcal
M}_c}.
\]
Composing the source and the target maps $\tilde{\alpha}$ and
$\tilde{\beta}$ of the cotangent groupoid $T^*G\gpd A^*G$ as defined
in Equation \eqref{eq:cotangent:groupoid} with $i_{{\mathcal
D}_{c}}^{*}: AG^{*} \to {\mathcal D}_{c}^{*}$, the dual map of the
canonical inclusion $i_{{\mathcal D}_{c}}: {\mathcal D}_{c} \to AG$,
we have the maps
\begin{equation}\label{eq:NH:source:target}
\begin{array}{l}
\tilde{\alpha}_{{\mathcal D}_c}=i^*_{{\mathcal D}_c}\circ
\tilde{\alpha}:T^*G\to {\mathcal D}_c^*\\ \tilde{\beta}_{{\mathcal
D}_c}=i^*_{{\mathcal D}_c}\circ \tilde{\beta}:T^*G\to {\mathcal
D}_c^*.
\end{array}
\end{equation}

If we apply the forward integrability algorithm to the implicit
difference equation $S_{(L,{\mathcal M}_c)}$ in $T^*G$, but changing
the maps $\tilde{\alpha}$ and $\tilde{\beta}$ for the new maps
$\tilde{\alpha}_{{\mathcal D}_c}$ and $\tilde{\beta}_{{\mathcal
D}_c}$ (see Remark \ref{rem:generalization}), then the first step of
the algorithm yields a close relation with the discrete
non-holonomic Euler-Lagrange equations. Indeed, if we have the
following result.

\begin{proposition}\label{forward-NHDEL}
Let $G\gpd M$ be a Lie groupoid and $(L,{\mathcal M}_c,{\mathcal
D}_c)$ be a discrete nonholonomic Lagrangian system on $G$. The
first step, $(S_{(L,{\mathcal M}_c)})^1_f$, to obtain the forward
integrable part of $S_{(L,{\mathcal M}_c)}=dL({\mathcal M}_c)$ for
the maps $\tilde{\alpha}_{{\mathcal D}_c}:T^*G\to {\mathcal D}^*_c$
and $\tilde{\beta}_{{\mathcal D}_c}:T^*G\to {\mathcal D}^*_c$,
defined in \eqref{eq:NH:source:target}, are the points $dL(g)\in
S_{(L,{\mathcal M}_c)}$ such that there exists an element $h\in
{\mathcal M}_c$, satisfying $(g, h)\in G_2$ and, in addition,
$(g,h)$ is a solution of the discrete non-holonomic Euler-Lagrange
equations, that is,
\[
d\left[L\circ l_{g}+L\circ r_{h}\circ i\right](\epsilon (x))_{|
({\mathcal D}_c)_x}=0,
\]
where $x=\beta (g)=\alpha (h)$.

In other words,
\[
\begin{array}{rcl}
(S_L)^1_f&=&\{ (g,h)\in G_2 \cap ({\mathcal M}_c\times {\mathcal M}_c)\, |\\
     && (g,h) \mbox{ is a solution of the discrete nonholonomic Euler-Lagrange eqns.}\}\\
&=&\{ (g,h)\in G_2\cap ({\mathcal M}_c\times {\mathcal M}_c)\, | \,
\tilde\beta_{{\mathcal D}_c} (dL(g))=\tilde\alpha_{{\mathcal
D}_c}(dL(h))\}
\end{array}
\]
\end{proposition}
\begin{proof}
The proof is analogous to that of Proposition \ref{forward-DEL},
taking into account the version of the discrete nonholonomic
Euler-Lagrange equations given by Equation
\eqref{eq:NHDEL:Legendre}.
\end{proof}
\subsection{Example: The discrete Chaplygin sleigh \cite{F, FZ}}

To illustrate the results contained in the previous section, we will describe
a discretization of the Chaplygin sleigh system (previously
considered in \cite{IMMM}).

The Chaplygin sleigh system describes the motion of a rigid body
sliding on a horizontal plane. The body is supported at three
points, two of which slide freely without friction while the third
is a knife edge, a constraint that allows no motion orthogonal to
this edge (see \cite{NF}).

The configuration space of this system is the group $SE(2)$ of Euclidean
motions of $\R^2$. An element $A \in SE(2)$ is represented by a
matrix
\[
A = \left(
\begin{array}{ccc}
\cos\theta &-\sin\theta & x\\
\sin\theta &\cos\theta & y\\
0&0&1
\end{array}
\right) \qquad\mbox{ with }\theta ,x,y\in \R .
\]
Thus, $(\theta, x, y)$ are local coordinates on $SE(2)$.

A basis of the Lie algebra $\mathfrak{se}(2)\cong \R ^3$ of $SE(2)$
is given by
\[
e=\left(
\begin{array}{ccc}
0&-1&0\\
1&0&0\\
0&0&0
\end{array}
\right),\qquad
e_1=\left(
\begin{array}{ccc}
0&0&1\\
0&0&0\\
0&0&0
\end{array}
\right) ,\qquad e_2=\left(
\begin{array}{ccc}
0&0&0\\
0&0&1\\
0&0&0
\end{array}
\right)
\]
and we have that
\[
[e,e_1]=e_2,\quad [e,e_2]=-e_1,\quad [e_1,e_2]=0.
\]
An element $\xi \in \mathfrak{se}(2)$ is of the form
\[
\xi =\omega \, e+v_1\, e_1+v_2\, e_2
\]
The exponential map $\exp: \mathfrak{se}(2)\cong \R^3\to
SE(2)$ of $SE(2)$ is given by
\[
\exp (\omega, v_1,v_2)=(\omega ,v_1\frac{\sin \omega}{\omega}+v_2 ( \frac{\cos \omega -1}{\omega } ),
-v_1 ( \frac{\cos \omega -1}{\omega} )+v_2  \frac{\sin \omega}{\omega }  ),\mbox{ if }\omega \neq 0,
\]
and
$\exp (0, v_1,v_2)=(0,v_1,v_2)$.
The restriction of the exponential map to the open subset $U=]-\pi ,\pi [\times \R^2\subseteq \R ^3
\cong \mathfrak{se}(2)$ is a diffeomorphism onto the open subset $\exp (U)$ of $SE(2)$.

A discretization of the Chaplygin sleigh may be constructed as follows
(see \cite{IMMM} for more details):

First of all, the discrete Lagrangian $L: SE(2)\longrightarrow \R$
is given by
\[
L(A)=\frac{1}{2} \hbox{Tr } ( A \J A^T)-\hbox{Tr
} ( A \J),
\]
where $\J$ is the matrix:
\[
\J=\left(
\begin{array}{ccc}
(J/2)+ma^2&mab&ma\\
mab&(J/2)+mb^2&mb\\
ma&mb&m
\end{array}
\right)
\]
(see \cite{FZ}). In terms of the coordinates $(\theta, x,y)$, the
discrete Lagrangian can be written as
\begin{eqnarray*}
L(\theta, x, y)&=&(max+mby-ma^2-mb^2-J)\cos\theta+m(ay-bx)\sin \theta\\&&+\frac{m}{2}\left((x-a)^2+(y-b)^2\right)+\frac{1}{2}(J-m)
\end{eqnarray*}
Second, the vector subspace ${\mathcal D}_c$ of the Lie
algebra $\mathfrak{se}(2)$ is given by
\[
{\mathcal D}_c=\hbox{ span } \{ e, e_1\}= \set{(\omega, v_1,
v_2)\in \mathfrak{se}(2)}{v_2=0}.
\]
Finally, the constraint submanifold ${\mathcal M}_c$ of
$SE(2)$ is ${\mathcal M}_c=\exp (U\cap {\mathcal D}_c)$, that is,
\begin{eqnarray*}
{\mathcal M}_c&=&\set{(\theta , x,y )\in SE(2)}{-\pi <\theta
<\pi,\, \theta \neq 0, (1-\cos \theta )x-y\sin \theta =0}\\ &&
\cup \set{(0,x,0)\in SE(2)}{x\in\R}.
\end{eqnarray*}

Let us calculate the set
$S_{(L,{\mathcal M}_c)}=dL({\mathcal M}_c)$.
We have that $(\theta, x, y; p_{\theta}, p_x, p_y)\in T^*SE(2)\equiv \R^6$ belongs to $S_{(L,{\mathcal M}_c)}$ if it verifies the following equations:
 \begin{eqnarray*}
 p_{\theta}&=&m(ay-bx)\cos\theta+(ma^2+mb^2+J-max-mby)\sin\theta, \\
p_x&=&ma\cos\theta-mb\sin\theta+m(x-a),\\
 p_y&=&mb\cos\theta+ma\sin\theta+m(y-b),\\
 \mbox{ with }(\theta, x, y)&\in& {\mathcal M}_c\; .
\end{eqnarray*}

Moreover, using \eqref{eq:cotangent:groupoid}, we deduce that
the maps $\tilde{\beta}:T^*SE(2)\to \mathfrak{se}(2)^\ast$ and
$\tilde{\alpha}:T^*SE(2)\to \mathfrak{se}(2)^\ast$ are just the pullbacks
by the left and right translations, respectively. Now, the mappings $\widetilde{\alpha}_{{\mathcal D}_c}: T^*G\to {\mathcal D}_c^*$ and $\widetilde{\beta}_{{\mathcal D}_c}: T^*G\to {\mathcal D}_c^*$ are:
\begin{eqnarray*}
\widetilde{\alpha}_{{\mathcal D}_c}(p_{\theta}d\theta +p_xdx+p_ydy)&=&(p_{\theta}-yp_x+xp_y)e^*+p_xe_1^*\\
\widetilde{\beta}_{{\mathcal D}_c}(p_{\theta}d\theta +p_xdx+p_ydy)&=&p_{\theta}e^*+(p_x\cos\theta+p_y\sin\theta)e_1^*
\end{eqnarray*}
where $p_{\theta}d\theta +p_xdx+p_ydy\in T^*_{(\theta, x, y)}G$ and $\{e^*, e^*_1\}$ is the dual basis of $\{e, e_1\}$ and therefore
$\hbox{span }\{e^*, e_1^*\}={\mathcal D}_c^*$.

Thus, the first step of the algorithm allows yields the set of points
$(\theta _k,x_k,y_k)\in {\mathcal M}_c$, with $k\in \{1,2\}$, such that
\[
\widetilde{\beta}_{{\mathcal D}_c} (dL(\theta _1,x_1,y_1))=\widetilde{\alpha}_{{\mathcal D}_c} (dL(\theta_2,x_2,y_2)),
\]
which are just the discrete
Euler-Poincar\'e-Suslov equations,
\begin{eqnarray*}
\left(\begin{array}{l}
-ma\cos\theta _1-mb\sin\theta _1+ma\\
+mx_1\cos\theta _1+my_1\sin\theta _1
\end{array}\right)
&\kern-5pt=&\kern-5pt
\left(\begin{array}{l}
mx_{2}+ma\cos\theta_{2}\\
-mb\sin\theta_{2}-ma \end{array}\right)
\\
\left(\begin{array}{l}
 m(ay_1-bx_1)\cos\theta_1-m(ax_1+by_1) \sin \theta_1
 \\
+(ma^2+mb^2+J)\sin \theta_1
\end{array}\right)
&\kern-5pt=&\kern-5pt
\left(\begin{array}{l} ma\,y_{2}-mb\,x_2\\
+(ma^2+mb^2+J)\sin \theta _2
\end{array}\right)
\end{eqnarray*}
where $(\theta _k,x_k,y_k)\in {\mathcal M}_c$, with $k\in \{1,2\}$.

\subsection{An approach to implicit discrete Hamiltonian  systems}

Along the paper, we have focused our attention to the case of discrete lagrangian mechanics including this
theory in the setting of discrete implicit systems. Of course, we can also adopt a dual point of view, that is,
the hamiltonian formalism. In the continuous setting, the hamiltonian formalism is specified given a hamiltonian
function $H$  on $A^*G$ equipped with its canonical linear Poisson bracket. The  dynamics is given by the
associated hamiltonian vector field $X_H$. In the sequel, we will show how to obtain a lagrangian submanifold of
the symplectic groupoid $T^*G\rightrightarrows A^*G$ directly from the flow of $X_H$. 

Assume that $G\rightrightarrows M$ is a Lie groupoid with corresponding algebroid $\tau: AG\rightarrow M$. Then,
we have already explained (see \eqref{eq:cotangent:groupoid}) that the cotangent bundle $T^*G$ is a groupoid
over $A^*G$. Moreover, there is a Poisson structure on $A^*G$ (induced from the Lie algebroid structure on $AG$)
such that the source map $\tilde\alpha:T^*G\to A^*G$ (resp. the target map $\tilde\beta:T^*G\to A^*G$) is a
Poisson (resp. anti-Poisson) map.

Let $H: A^*G\longrightarrow \R$ be a Hamiltonian function. Then, since $A^*G$ is a Poisson manifold, we have the
hamiltonian vector field $X_H$. On the other hand, consider the pull-back of the 1-form $dH\in \Omega^1(A^*G)$
by the source map $\tilde{\alpha}$. Using the canonical symplectic structure on $T^*G$, we have the
corresponding Hamiltonian vector field $X_{\tilde\alpha^*H}$. It is easy to see that both vector fields are
$\tilde\alpha$-related. Indeed, if $\Pi _{A^*G}$ (resp. $\Pi_G$) denotes the Poisson structure on $A^*G$ (resp.
on $T^*G$), since $\tilde{\alpha}_* (\Pi_G)=\Pi _{A^*G}$,
\[
X_H=\Pi _{A^*G}^\sharp (dH)=(\tilde\alpha_*\Pi _G )^\sharp (dH)= \tilde\alpha_* ((\Pi _G )^\sharp (d
\tilde\alpha ^*H))=\tilde\alpha_* (X_{\tilde\alpha^*H}),
\]
where $\Pi ^\sharp (dF)$ denotes contraction by $dF$, i.e., $\Pi ^\sharp (dF)(\gamma)=\Pi (dF, \gamma)$ for any
1-form $\gamma$. Thus, we have that the following diagram is commutative

\hspace{-2cm}\begin{picture}(375,60)(40,50)
\put(200,20){\makebox(0,0){$A^*G$}} \put(305,20){\makebox(0,0){$A^*G$}} \put(225,20){\vector(1,0){55}}
\put(200,80){\makebox(0,0){$T^*G$}} \put(200,70){\vector(0,-1){40}} \put(300,70){\vector(0,-1){40}}
\put(225,80){\vector(1,0){55}} \put(300,80){\makebox(0,0){$T^*G$}}

\put(250,85){$\varphi ^t _{H}$} \put(250,25){$\psi ^t _{H}$}

\put(190,50){\makebox(0,0){$\tilde\alpha$}} \put(310,50){\makebox(0,0){$\tilde\alpha$}}
\end{picture}

\vspace{1.5cm} \noindent where $\varphi ^t _{H}$ (resp., $\psi ^t_{H}$) is the flow of $X_{\tilde\alpha^*H}$
(resp., $X_H$). In addition, since $\varphi ^t _{H}:T^*G\to T^*G$ is a symplectomorphism, $\tilde{i}:T^*G\to T^*G$ is an anti-symplectomorphism and
$\tilde{\epsilon}:A^*G\to T^*G$ is a Lagrangian embedding, we deduce that the submanifold $S^t_H$ defined by
\[
S^t_{H}=\left\{ \nu _g\in T^*G\; |\; \nu _g = \tilde{i}(\varphi^t_{H}(\tilde\epsilon(x))), \forall x\in A^*G\right\}\subset
T^*G,
\]
is a Lagrangian submanifold of $T^*G$.
\begin{remark}
We must note that the previous discussion is a particular instance
of the construction of Lagrangian bisections on symplectic Lie
groupoids (see \cite{CDW}).
\end{remark}

\begin{example}{\rm 
Let $H:T^*Q\to \R$ be a Hamiltonian function on $T^*Q$. In this case, using canonical coordinates $(q,p)$, the
Poisson structure on $T^*Q$ is $\Pi _Q=-\frac{\partial}{\partial q}\wedge \frac{\partial}{\partial p}$ and the
Hamiltonian vector field is given by
\[
X_H =\frac{\partial H}{\partial p}\frac{\partial}{\partial q}-\frac{\partial H}{\partial q}
\frac{\partial}{\partial p}.
\]
Now, let $G=Q\times Q$ be the pair groupoid. Then, on $T^*(Q\times Q)$ the source map $\tilde\alpha :T^*(Q\times
Q)\to T^*Q$ is just $(q_0,q_1,p_0,p_1)\mapsto (q_0,-p_0)$. As a consequence, since $\Pi_{Q\times
Q}=\frac{\partial }{\partial q_0}\wedge\frac{\partial}{\partial p_0}+\frac{\partial }{\partial q_1}\wedge
\frac{\partial}{\partial p_1}$, we have
\[
\begin{array}{rcl}
d\tilde\alpha^* H&=&\displaystyle \frac{\partial H}{\partial q_0}dq_0-\frac{\partial H}{\partial p_0}
dp_0  \\[8pt]
X_{\tilde\alpha^* H}&=&\displaystyle \frac{\partial H}{\partial p_0} \frac{\partial}{\partial
q_0}+\frac{\partial H}{\partial q_0}\frac{\partial}{\partial p_0}
\end{array}
\]
Finally, since $\tilde\epsilon :A^*G\to T^*G$, $(q,p)\mapsto (q,q,-p,p)$, we have that
\[
S^t_H=\left\{ (q,(\psi _1)_H ^t (q,p),-p,(\psi _2)_H ^t (q,p))\,|\, (q,p)\in T^*Q\right\},
\]
where $(q,p)\mapsto \psi ^t_H (q,p)=\Big( (\psi _1)_H ^t (q,p),(\psi _2)_H ^t (q,p)\Big)$ is the flow of $X_H$.}
\end{example}

\begin{remark}
Let $L:TQ\to {\Bbb R}$ be a regular continuous Lagrangian on the manifold $Q$.  \emph{The exact discrete Lagrangian  for a small  time-step $h>0$}  (\cite{MW}) is the discrete Lagrangian $L_d^h:Q\times Q\to {\Bbb R}$ given by 
$$L_d^h(q_0,q_1)=\int_0^hL(\sigma_{01}(t), \dot{\sigma}_{01}(t))dt$$
where $\sigma_{01}(t)$ is the unique solution of the Euler-Lagrange equations for $L$ which satisfies the boundary conditions $\sigma_{01}(0)=q_0$ and $\sigma_{01}(h)=q_1.$ 

The Legendre transformations, ${\Bbb F}^\pm_{L_d^h}$ and ${\Bbb F}L$, for $L_d^h$ and $L$ respectively, are related as follows (see \cite{MW})
$$
\begin{array}{rcl}
{\Bbb F}^+{L_d^h}(q_0,q_1)&=&{\Bbb F}L(\sigma_{01}(h),\dot{\sigma}_{01}(h))\\[5pt]
{\Bbb F}^-{L_d^h}(q_0,q_1)&=&{\Bbb F}L(\sigma_{01}(0),\dot{\sigma}_{01}(0))
\end{array}
$$
Therefore, if $L$ is (hyper-)regular then $L_d^h$ is (hyper-)regular. 

Suppose that $L$ is hyperregular. Denote by $H:T^*Q\to {\Bbb R}$ the corresponding Hamiltonian function, i.e.
$$H=E_L\circ {\Bbb F}L^{-1}$$
$E_L$ being the Lagrangian energy associated with $L.$

In \cite{MW} it is proved that for a small time-step $h$, the following diagrams relate the flow $\psi_H^h: T^*Q\to T^*Q$ at time $h$   of the Hamiltonian vector field $X_H\in {\frak X}(T^*Q)$  and  the discrete Lagrangian evolution operator $\Upsilon_{L_d^h}:Q\times Q\to Q\times Q$ associated with $L_d^h$ defined  as in (\ref{operador})

\vspace{0cm}

\hspace{-4cm}\begin{picture}(375,60)(40,50)
\put(200,20){\makebox(0,0){$T^*Q$}} \put(305,20){\makebox(0,0){$T^*Q$}} \put(225,20){\vector(1,0){55}}
\put(200,80){\makebox(0,0){$Q\times Q$}} \put(200,70){\vector(0,-1){40}} \put(300,70){\vector(0,-1){40}}
\put(225,80){\vector(1,0){55}} \put(300,80){\makebox(0,0){$Q\times Q$}}

\put(250,88){$\Upsilon_{L_d^h}$} \put(250,25){$\psi ^t _{H}$}

\put(190,50){\makebox(0,0){${\Bbb F}_{L_d^h}^+$}} \put(312,50){\makebox(0,0){${\Bbb F}_{L_d^h}^+$}}
\end{picture}

\vspace{-2.3cm}

\hspace{2cm}
\begin{picture}(375,60)(40,50)
\put(200,20){\makebox(0,0){$T^*Q$}} \put(305,20){\makebox(0,0){$T^*Q$}} \put(225,20){\vector(1,0){55}}
\put(200,80){\makebox(0,0){$Q\times Q$}} \put(200,70){\vector(0,-1){40}} \put(300,70){\vector(0,-1){40}}
\put(225,80){\vector(1,0){55}} \put(300,80){\makebox(0,0){$Q\times Q$}}

\put(250,88){$\Upsilon_{L_d^h}$} \put(250,25){$\psi ^t _{H}$}

\put(190,50){\makebox(0,0){${\Bbb F}_{L_d^h}^-$}} \put(312,50){\makebox(0,0){${\Bbb F}_{L_d^h}^-$}}
\end{picture}

\vspace{1.3cm}

A direct computation, using this result,  shows that 
$$\Psi(S_{L^h_d})=(\Psi\circ d{L_d^h})(Q\times Q)={\rm graph} \; \psi_H^h= \Psi(S_H^h)$$
Here $\Psi:T^*(Q\times Q)\to T^*(Q\times Q)$ is the isomorphism described in (\ref{Psi}). Thus, $S_{L_d^h}=S_H^h$. 

It would be interesting to extend this result for a general regular Lagrangian $L:AG\to {\Bbb R}$ on the associated Lie algebroid $AG$  of a Lie groupoid $G \rightrightarrows  M$. Previously, its is necessary to describe the exact discrete Lagrangian $L_d^h:G\to {\Bbb R}$ induced by $L$ (see \cite{MMM2} for more details). 

The results of this section motivate that  we may  introduce, in the general context of a Lie groupoid $G$, the notion of an implicit discrete Hamiltonian system as a Lagrangian submanifold of $T^*G.$

\end{remark}

\section{Conclusion}
We have developed an algorithm which permits to extract the integrable part of an implicit
difference equation. This algorithm produces a sequence of submanifolds (or generally
subsets) which encodes where
there exist well defined solutions of the discrete dynamics. As an application, our method
allows us to easily
analyze the case of discrete (nonholonomic) Lagrangian systems from an implicit point of
view, in particular the situation of singular lagrangians, and we will expect that it will be an useful tool for discrete optimal control theory. We
will discuss this topic in a future paper.

Finally, in our opinion,  the developments about implicit systems defined by Hamiltonian functions open
doors to future research in the theory of generating functions for Poisson morphisms (using standard symplectic
techniques), the theory of discrete hamiltonian systems in a Poisson context \cite{AhPe}, the Dirac theory of
constraints in the discrete setting, a discrete Hamilton-Jacobi theory...  among other research topics.

\appendix
\section{Discrete Lagrangian Mechanics on Lie groupoids}\label{A}

A \emph{discrete Lagrangian system} consists of a Lie groupoid
${G}\rightrightarrows M$ (the \emph{discrete space}) and a
\emph{discrete Lagrangian} $L: {G} \to \R$.

For $g\in {G}$ fixed, we consider the set of \emph{admissible
sequences}:
\[
{\mathcal C}^N_{g}\kern-2.1pt=\kern-2pt\set{ \gamma_G:[1,N]\cap\Z\to G} { \gamma_G \mbox{ is an admissible
sequence and } \gamma_G(1) .. \gamma_G(N)=g }.
\]
We may  identify the tangent space to ${\mathcal C}^N_g$ at
$\gamma_G$ with
\[
T_{\gamma_G}{\mathcal C}^N_g\equiv\left\{(v_1, v_2, \ldots, v_{N-1})\, |\, v_k\in A_{x_k}G \hbox{ and }
x_k=\beta(\gamma_G(k)), 1\leq k\leq N-1\right\}.
\]
$(v_1, v_2, \ldots, v_{N-1})$ is called an {\em infinitesimal
variation} of $\gamma_G$.

Now, we define the \emph{discrete action sum} associated to the
discrete Lagrangian $L: {G}\to \R$ by
\[
{\mathcal S} L (\gamma _G)= \sum_{k=1}^{N} L(\gamma_G(k)).
\]
Hamilton's principle requires that this discrete action sum be
stationary with respect to all the infinitesimal variations. This
requirement gives the following alternative expressions for the
\emph{discrete Euler-Lagrange equations} (see \cite{MMM}):
\begin{eqnarray}
  \lvec{X}\big({g_k})(L)-\rvec{X}\big({g_{k+1}})(L)=0,
  \label{discretee}
\end{eqnarray}
for all sections $X$ of $\tau: A{G}\to M$. Alternatively, we may rewrite
the discrete Euler-Lagrange equations as
\[
d\left[L\circ l_{g_k}+L\circ r_{g_{k+1}}\circ i\right](\epsilon (x_k))_{\big|A_{x_k}G}=0,
\]
where $\beta(g_k)=\alpha(g_{k+1})=x_k$.

Thus, we may define the \emph{discrete Euler-Lagrange operator} $D_{\hbox{\footnotesize DEL}}L: {G}_2\to A^*{G}$
from $G_2$ to $A^*G$, the dual of $A{G}$. This operator is given by
\[
D_{\hbox{\footnotesize DEL}}{L}(g, h)=  d\left[L\circ l_{g}+L\circ r_{h}\circ i\right](\epsilon
(x))_{\big|A_{x}G}
\]
with $\beta(g)=\alpha(h)=x$.

Let $\Upsilon: {G}\to {G}$ be a smooth map such that:
\begin{enumerate}
\item[-] $\hbox{graph}(\Upsilon)\subseteq {G}_2$, that is, $(g,
\Upsilon(g))\in
  {G}_2$, for all $g\in {G}$ ($\Upsilon$ is a \emph{second order operator}).
\item[-] $(g, \Upsilon(g))$ is a solution of the discrete
Euler-Lagrange
  equations, for all $g\in {G}$, that is,
  $(D_{\hbox{DEL}}L)(g,\Upsilon(g))=0,$ for all $g\in {G}.$
\end{enumerate}
In such a case
\begin{equation}\label{5.22}
  \lvec{X}(g)(L)-\rvec{X}(\Upsilon(g))(L)=0,
\end{equation}
for every section $X$ of $A{G}$ and every $g\in {G}.$ The map
$\Upsilon: {G}\to {G}$ is called a \emph{discrete flow} or a
\emph{discrete Lagrangian evolution operator for $L$}.

Given a Lagrangian $L: {G}\to \R$ we define the \emph{discrete
Legendre transformations} $\F^{-}L: {G}\to A^*{G}$ and $\F^{+}L:
{G}\to A^*{G}$ by
\begin{equation}\label{eq:Leg:transform}
\begin{array}{rcl}
  (\F^{-}L)(h)(v_{\epsilon(\alpha(h))})&= &-v_{\epsilon(\alpha(h))}(L\circ
  r_h\circ i) \\ &=&\tilde{\alpha}(dL(h)) (v_{\epsilon(\alpha(h))}),\;\;\; \mbox{ for }
  v_{\epsilon(\alpha(h))}\in A_{\alpha(h)}{G},\\[7pt]
  (\F^{+}L)(g)(v_{\epsilon(\beta(g))})&=& v_{\epsilon(\beta(g))}(L\circ
  l_g)\\&=&\tilde{\beta}(dL(g)) (v_{\epsilon(\beta(g))}) , \mbox{ for } v_{\epsilon(\beta(g))}\in
  A_{\beta(g)}{G}.
\end{array}
\end{equation}
A discrete Lagrangian $L:{G}\to \R$ is said to be \emph{regular} if and only if the Legendre transformation
$\F^+L$  is a local diffeomorphism (equivalently, if and only if the Legendre transformation $\F^-L$ is a local
diffeomorphism). In this case, if $(g_0,h_0)\in G\times G$ is a solution of the discrete Euler-Lagrange
equations for $L$ then, one may prove (see \cite{MMM}) that there exist two open subsets $U_{0}$ and $V_{0}$ of
$G$, with $g_{0} \in U_{0}$ and $h_{0} \in V_{0}$, and there exists a (local) discrete Lagrangian evolution
operator $\Upsilon_{L}: U_{0} \to V_{0}$ such that:
\begin{enumerate}
\item $\Upsilon_{L}(g_{0}) = h_{0}$, \item $\Upsilon_{L}$ is a diffeomorphism and \item $\Upsilon_{L}$ is
unique, that is, if $U'_{0}$ is an open subset of $G$, with $g_{0} \in U_{0}',$ and $\Upsilon'_{L}: U'_{0} \to
G$ is a (local) discrete Lagrangian evolution operator then
\[
\Upsilon_{L |U_{0}\cap U_{0}'} = \Upsilon'_{L |U_{0}\cap U_{0}'}.
\]
\end{enumerate}
Moreover, if $\F^{+}L$ and $\F^{-}L$ are global diffeomorphisms (that is, $L$ is \emph{hyperregular}) then
$\Upsilon _L=(\F^{-}L)^{-1}\circ \F^{+}L$.

If $L: {G}\to \R$ is a hyperregular Lagrangian function, then pushing forward to $A^*{G}$ with the discrete
Legendre transformations, we obtain the \emph{discrete  Hamiltonian evolution operator}, $\tilde{\Upsilon}_{L}:
A^*{G}\to A^*{G}$ which is given by
\begin{equation}\label{dheo}
\tilde{\Upsilon}_{L}=\F^{\pm}L\circ \Upsilon_{L}\circ
(\F^{\pm}L)^{-1}\;=\F^{+}L\circ (\F^{-}L)^{-1}.
\end{equation}

\section{Discrete nonholonomic Lagrangian systems}

\emph{ A discrete nonholonomic Lagrangian system } on a Lie groupoid $G\rightrightarrows M$ is a 
\emph{Lagrangian discrete} $L:G\to {\Bbb R}$, a vector subbundle ${\mathcal D}_c$ (\emph{the constraint distribution}) of the Lie algebroid $AG$ of $G$ and a \emph{ discrete constraint embedded  submanifold } 
${\mathcal M}_c$ of $G$ such that $\dim {\mathcal M}_c=\dim {\mathcal D}_c.$

Let $(L,{\mathcal M}_c,{\mathcal D}_c)$ be a discrete nonholonomic
Lagrangian system. The \emph{ discrete nonholonomic Euler-Lagrange
equations for the system $(L_d,{\mathcal M}_c,{\mathcal D}_c)$} are
given by
\[
d(L_d\circ l_g + L_d\circ r_h\circ i)(\epsilon (x))_{|({\mathcal
D}_c)_x}=0,
\]
for $(g,h)\in G_2\cap ({\mathcal M}_c\times {\mathcal M}_c)$, with
$\beta(g)=\alpha(h)=x$ (for more details on this section, see
\cite{IMMM}).

For a discrete nonholonomic Lagrangian system, we can define the
\emph{discrete nonholonomic Legendre transformations}
\[
\F^-(L, {\mathcal M}_{c}, {\mathcal D}_{c}): {\mathcal M}_{c} \to
{\mathcal D}_{c}^{*} \; \; \mbox{ and } \; \; \F^+(L, {\mathcal
M}_{c}, {\mathcal D}_{c}): {\mathcal M}_{c} \to {\mathcal D}_{c}^{*}
\]
as follows:
\begin{eqnarray}
\F^{-}(L, {\mathcal M}_{c}, {\mathcal
D}_{c})(h)(v_{\epsilon(\alpha(h))})&\kern-5pt=
&\kern-5pt-v_{\epsilon(\alpha(h))}(L \circ r_{h} \circ i), \nonumber
\\\label{F-LMD}
&\kern-5pt= &\kern-5pt \tilde{\alpha}(dL(h))
(v_{\epsilon(\alpha(h))}), \mbox{ for } v_{\epsilon(\alpha(h))} \in
{\mathcal D}_{c}(\alpha(h)),
\\
 \F^{+}(L, {\mathcal M}_{c}, {\mathcal
D}_{c})(g)(v_{\epsilon(\beta(g))})&\kern-5pt= &\kern-5pt
v_{\epsilon(\beta(g))}(L \circ l_{g}), \nonumber
\\\label{F+LMD} &\kern-5pt=
& \kern-5pt\tilde{\beta}(dL(g)) (v_{\epsilon(\beta(g))}), \mbox{ for
} v_{\epsilon(\beta(g))} \in {\mathcal D}_{c}(\beta(g)).
\end{eqnarray}
If $\F^{-}L_{d}: G \to AG^{*}$ and $\F^{+}L_{d}: G \to AG^{*}$ are
the standard discrete Legendre transformations associated with the
Lagrangian function $L$ defined in Equation \eqref{eq:Leg:transform}
and $i_{{\mathcal D}_{c}}^{*}: AG^{*} \to {\mathcal D}_{c}^{*}$ is
the dual map of the canonical inclusion $i_{{\mathcal D}_{c}}:
{\mathcal D}_{c} \to AG$ then
\begin{equation}\label{Rel}
\F^{-}(L_{d}, {\mathcal M}_{c}, {\mathcal D}_{c}) = i_{{\mathcal
D}_{c}}^{*} \circ \F^{-}L_{d} \circ i_{{\mathcal M}_{c}}, \; \;
\F^{+}(L_{d}, {\mathcal M}_{c}, {\mathcal D}_{c}) = i_{{\mathcal
D}_{c}}^{*} \circ \F^{+}L_{d} \circ i_{{\mathcal M}_{c}}.
\end{equation}
Using Equations \eqref{F-LMD} and \eqref{F+LMD}, the discrete
nonholonomic Euler-Lagrange equations are equivalent to
\begin{equation}\label{eq:NHDEL:Legendre}
\F^{-}(L_{d}, {\mathcal M}_{c}, {\mathcal D}_{c})(h)= \F^{+}(L_{d},
{\mathcal M}_{c}, {\mathcal D}_{c})(g).
\end{equation}



\begin{thebibliography}{99}

\bibitem{AhPe} C. D. Ahlbrandt, A. P. Peterson: {\sl Discrete Hamiltonian Systems, Difference equations, Contunued Fractions and Riccati Equations,} Kluwer Texts in the Mathematical Sciences, Kluwer, Netherlands, 1996.

\bibitem{Ca} J. F. Cari\~{n}ena: Theory of singular
lagrangians, {\sl Fortschr. Phys. } {\bf 38} (1990) 641--679.

\bibitem{CeEt} H. Cendra, M.  Etchechoury: Desingularization of implicit analytic differential equations. {\sl J. Phys. A} {\bf 39} (35) (2006),  10975-11001.


\bibitem{CMR}
H. Cendra, J. E. Marsden, T.S. Ratiu: {\sl Lagrangian reduction by stages}, {\rm Mem. Amer. Soc.} {\bf 152}
(2001) no. 722.

\bibitem{CDW}
A. Coste, P. Dazord, A. Weinstein: Grupo\"\i des symplectiques, {\sl Pub. D\'{e}p. Math. Lyon}, {\bf 2/A}
(1987), 1--62.

\bibitem{CL} M. Crainic, R.L. Fernandes: Integrability of Lie
brackets. {\sl Ann. of Math.}  {\bf 157} (2) (2003), 575--620.

\bibitem{F} Y.N. Fedorov: A Discretization of
the Nonholonomic Chaplygin Sphere Problem. {\sl SIGMA} {\bf 3} (2007), 044, 15 pages.

\bibitem{FZ} Y.N. Fedorov, D.V. Zenkov: Discrete
nonholonomic LL systems on Lie groups, {\sl Nonlinearity} {\bf 18}
(2005), 2211--2241.


\bibitem{Hair}
E. Hairer, C. G. Wanner: {\sl Geometric Numerical Integration, Structure-Preserving Algorithms for Ordinary
Differential Equations}, Springer Series in Computational Mathematics {\bf 31}, Springer-Verlag, Berlin
Heidelberg, 2002.

\bibitem{IMMM}
D. Iglesias, J.C. Marrero, D. Mart\'{\i}n de Diego, E.
Mart{\'\i}nez: Discrete nonholonomic Lagrangian systems on Lie
groupoids. {\sl J. Nonlinear Sci}.  {\bf 18} (3)  (2008),  221--276.

\bibitem{IMMS}
D. Iglesias, J.C. Marrero, D. Mart\'{\i}n de Diego, D. Sosa: Singular Lagrangian systems and variational
constrained mechanics on Lie algebroids, {\sl Dynamical Systems}, {\bf 23} No. 3 (2008), 351--397.


\bibitem{LMM}
M. de Le\'{o}n, J.C. Marrero, E. Mart\'{\i}nez: Lagrangian submanifolds and dynamics on Lie algebroids, {\sl J.
Phys. A: Math. Gen.} \textbf{38} (2005) R241--R308.

\bibitem{Ma} K. Mackenzie:
General theory of Lie groupoids and Lie algebroids. {\sl London Mathematical Society Lecture Note Series}, {\bf
213}, Cambridge University Press, Cambridge, 2005.

\bibitem{MMT} G. Marmo, G. Mendella, W. M. Tulczyjew:
Symmetries and constants of the motion for dynamics in implicit form, {\sl Ann. Inst. Henri Poincar\'e}, {\bf
57} 2, (1992), 147--166.

\bibitem{MMT1} G. Marmo, G. Mendella, W. M. Tulczyjew: Integrability
of implicit differential equations. {\sl  J. Phys. A: Math. Gen.} {\bf 30} (1), (1995), 149--163.

\bibitem{MMT2} G. Marmo, G. Mendella, W. M. Tulczyjew: Constrained Hamiltonian
systems as implicit differential equations. {\sl  J. Phys. A: Math. Gen.} {\bf 30}  (1), (1997) 277--293.


\bibitem{MW} J.E. Marsden, M. West: Discrete mechanics and variational integrators, {\sl Acta Numerica,} (2001) 357--514. 

\bibitem{MMM}
J.C. Marrero, D. Mart\'{\i}n de Diego, E. Mart\'{\i}nez: Discrete Lagrangian and Hamiltonian Mechanics on Lie
groupoids, {\sl Nonlinearity} {\bf 19} (2006), 1313--1348. Corrigendum: {\sl Nonlinearity} {\bf 19} (2006),
3003--3004.

\bibitem{MMM2}
J.C. Marrero, D. Mart\'{\i}n de Diego, E. Mart\'{\i}nez: The exact discrete Lagrangian function on Lie groupoids and some applications,  work in progress.

\bibitem{MMS}
J.C. Marrero, D. Mart\'{\i}n de Diego, A. Stern: Lagrangian submanifolds and discrete constrained mechanics on
Lie groupoids. {\sl Preprint} 2010.


\bibitem{NF}
 J. Neimark, N. Fufaev: {\sl Dynamics on Nonholonomic systems},
Translation of Mathematics Monographs, {\bf 33}, AMS,
Providence, RI, 1972.

\bibitem{PeLo} L. Petzold, P. L\"{o}tstedt: Numerical solution of nonlinear differential equations with algebraic constraints. II. Practical implications. {\sl SIAM J. Sci. Statist. Comput.} {\bf 7} (3) (1986),  720-733.

    \bibitem{RaRh} P. J. Rabier, W.C. Rheinboldt: Finite Diffreence Methods for Time Dependent, Linear Differential Algebraic equations, {\sl Appl. Math. Lett.} {\bf 7} (2) (1994), 29--34.


\bibitem{Sanz}
J. M. Sanz-Serna,  M. P. Calvo: {\sl Numerical Hamiltonian Problems}, Chapman\& Hall, London 1994

\bibitem{Stern}
A. Stern: Discrete {H}amilton--{P}ontryagin mechanics and generating functions on {L}ie groupoids, {\sl J. Symplectic Geom.} {\bf 8} n. 2 (2010), 225-238.


\bibitem{Tul1} W. Tulczyjew: Les sous-vari\'{e}t\'{e}s lagrangiennes et la
dynamique hamiltonienne, {\sl C.R. Acad. Sci. Paris} {\bf 283} (1976), 15-18.

\bibitem{Tul2} W. Tulczyjew: Les sous-vari\'{e}t\'{e}s lagrangiennes et la
dynamique lagrangienne, {\sl C.R. Acad. Sci. Paris} {\bf 283} (1976), 675-678.

\bibitem{weinstein:coisotropic}
A. Weinstein: Coisotropic calculus and Poisson groupoids, {\sl J. Math. Soc. Japan} {\bf 40} (1988), 705-727.

\bibitem{weinstein}
A. Weinstein: Lagrangian Mechanics and groupoids, {\sl Fields Inst. Comm.} {\bf 7} (1996), 207-231.




\end{thebibliography}
\end{document}